\documentclass{article}

\usepackage{graphicx}

\usepackage{amsmath,amsfonts,amssymb,latexsym,epsfig}
\usepackage{mathrsfs}
\usepackage{verbatim}
\usepackage{latexsym}
\usepackage{amsthm}
\usepackage{amssymb}
\usepackage{subfig}
\usepackage{graphics}
\usepackage{amsbsy}
\usepackage{fullpage}
\usepackage{enumerate}
\usepackage{times}
\usepackage{multirow}
\usepackage{soul}
\usepackage[normalem]{ulem}

\newcommand{\knot}{k_0}
\newcommand{\ktnot}{\widetilde{k}_0}
\newcommand{\kh}{k_h}
\newcommand{\kth}{\widetilde{k}_h}
\newcommand{\xnot}{x_0}
\newcommand{\xtnot}{\widetilde{x}_0}
\newcommand{\xh}{x_h}
\newcommand{\xth}{\widetilde{x}_h}
\newcommand{\xone}{x_1}
\newcommand{\xtone}{\widetilde{x}_1}
\newcommand{\dnot}{\delta_0}
\newcommand{\deltah}{\delta_h}
\newcommand{\done}{\delta_1}
\newcommand{\Dnot}{\Delta_0}
\newcommand{\Dh}{\Delta_h}
\newcommand{\Lp}{{L^\prime}}
\newcommand{\R}{\mathbb R}

\newtheorem{theorem}{Theorem}[section]

\newtheorem{remark}[theorem]{Remark}

\numberwithin{equation}{section}
\usepackage{amscd}
\usepackage{tikz}
\usepackage{xcolor}
\usepackage[english]{babel}
\usepackage[latin1]{inputenc}
\usepackage{times}
\usepackage[T1]{fontenc}
\usepackage{graphicx}
\usepackage{amsmath,amssymb}
\usepackage{slidesec}
\usepackage{verbatim}
\usepackage{url}
\usepackage{epsf}
\usepackage{amsmath}
\usepackage{graphics}
\usepackage{amsfonts}
\usepackage{amsbsy}
\usepackage{lscape}
\usepackage{enumerate}
\usepackage{amsthm}
\usepackage{amssymb}
\usepackage{exscale}
\usepackage{algorithm}
\usepackage{algorithmic}

\newcommand{\modkz}[1]{\textcolor{black}{#1}\index {#1}}
\newcommand{\modkzz}[1]{\textcolor{black}{#1}\index {#1}}
\newcommand{\modk}[1]{\textcolor{black}{#1}\index {#1}}

\begin{document}
\title{Contractivity of Runge-Kutta methods for convex gradient systems}
\author{J. M. Sanz Serna $^{1}$ \and Konstantinos C. Zygalakis$^{2}$}

\maketitle

\begin{abstract}
We consider the application of Runge-Kutta (RK) methods to gradient systems \((d/dt)x = -\nabla
V(x)\), where, as in many optimization problems, \(V\) is convex and \(\nabla V\) (globally)
Lipschitz-continuous with Lipschitz constant \(L\). Solutions of this system behave contractively, i.e.\ the
Euclidean
distance between two solutions \(x(t)\) and \(\widetilde{x}(t)\) is a nonincreasing function of \(t\). It is
then of interest to investigate whether a similar contraction takes place, at least for suitably small step
sizes \(h\), for the discrete solution.  \modkz{Dahlquist and Jeltsch 
results'  imply  that (1) there are explicit RK schemes that behave contractively whenever \(Lh\) is below a scheme-dependent constant and (2) Euler's rule is optimal in this regard. We prove however, by explicit construction of a convex potential using ideas from robust control theory,  that there exists RK schemes that fail to behave contractively for any choice of the time-step  \(h\).}
\end{abstract}
\footnotetext[1]{Departamento de Matem\'aticas, 
	Universidad Carlos III de Madrid, 
 	Legan\'es (Madrid), Spain}
\footnotetext[2]{School of Mathematics, University of Edinburgh, Edinburgh, Scotland}

\section{Introduction}
Systems of differential equations
\begin{equation}\label{eq:generalode}
\frac{d}{dt} x = F(x),
\end{equation}
with the gradient structure
\begin{equation}\label{eq:gradient}
\frac{d}{dt} x = -\nabla V(x),
\end{equation}
arise in many applications
and, accordingly,  have  attracted the interest of numerical analysts for a long time, see e.g.
\cite{HS94,HL14} among many others.  Here \(V\) is a continuously differentiable real function defined in \(\R^d\); in optimization applications \(V\) is the objective function and in Physics problems corresponds to a potential.
Since \((d/dt)V(x(t))\leq 0\),  \(V\) decreases along solutions. Furthermore, if
\(\lim_{t\rightarrow \infty} x(t) = x^\star\), then \(x^\star\) is a stationary point of \(V\), i.e., \(\nabla
V(x^\star) = 0\). These facts explain the well-known connections between  numerical integrators  for
\eqref{eq:gradient} and algorithms for the minimization of \(V\). The simplest example is provided by the
Euler integrator, that gives rise to the gradient descent optimization algorithm \cite{N14}. In the case where
\(\nabla V\) possesses a global Lipschitz constant \(L>0\) and \eqref{eq:gradient} is integrated with an
\emph{arbitrary} Runge-Kutta (RK) method, Humphries and Stuart \cite{HS94} showed that
 the value of \(V\) decreases along the computed solution, i.e.\ \(V(x_{n+1})\leq
V(x_n)\), for  positive stepsizes \(h\) with \(h\leq h_0\), where \(h_0>0\) only depends on \(L\) and on the
RK scheme.

In view of the important role that \emph{convex} objective functions play in optimization theory, see e.g.
\cite[Section 2.1.2]{N14}, it is  certainly of  interest to study numerical integrators for
\eqref{eq:gradient} in the specific case where \(V\) is convex, i.e.,
\begin{equation}\label{eq:onesided}
\forall x,y,\qquad \langle \nabla V(x)-\nabla V(y), x-y\rangle \geq 0
\end{equation}
(\(\langle \cdot,\cdot\rangle\) and \(\|\cdot\|\) stand throughout for the Euclidean inner product and norm in
\(\R^d\)). After recalling (see \cite[Section IV.2]{HaW96} or \cite[Definition 112A]{B16})  that a system of the general form \eqref{eq:generalode} is said to have
\emph{one-sided Lipschitz constant} \(\nu\) if
\begin{equation}\label{eq:nu}
\forall x,y,\qquad \langle F(x)-F(y), x-y\rangle \leq \nu \| x-y\|^2,
\end{equation}
we conclude that, for \modkz{convex} gradient systems \eqref{eq:gradient}, \(\nu = 0\).
 It follows that, for any two solutions \(x(t)\),
\(\widetilde x(t)\) of a gradient system, we have the \emph{contractivity estimate}
\begin{equation}\label{eq:contractivity}
\forall t\geq 0,\qquad \|\widetilde  x(t)-x(t)\|\leq \|\widetilde x(0)-x(0)\|,
\end{equation}
and in particular for any solution \(x(t)\) and any stationary point \(x^\star\) (which by convexity will
automatically be a minimizer)
\[
\forall t\geq 0,\qquad \|x(t)-x^\star\|\leq \|x(0)-x^\star\|.
\]

The study of linear multistep methods that, when applied to systems  of the general form \eqref{eq:generalode}
with one-sided Lipschitz constant \(\nu=0\), mimic the contractive behaviour in \eqref{eq:contractivity} began with the pioneering work of Dahlquist \cite{GD76}. The corresponding results in the Runge-Kutta (RK) field followed immediately \cite{JB75}.
Those developments gave rise to the notions of algebraic stability/B-stability of RK methods (see \cite[Section IV.12]{HaW96}, \cite[Section 357]{B16} and the monograph \cite{DK}) and
G-stability of  multistep methods (\cite[Section V.6]{HaW96} or \cite[Section 45]{B16}). These notions
 extend the concepts of A-stability \cite{GD63}
 to a nonlinear setting. Of course, algebraically stable/B-stable RK schemes and G-stable
multistep methods have to be \emph{implicit} and therefore are not well suited to be the basis of optimization
algorithms for large problems.

In this article we focus on the application of RK methods to gradient systems \eqref{eq:gradient} where  \(V\) is
convex and \(\nabla V\) is Lipschitz continuous with Lipschitz constant \(L\), i.e.
\[
\forall x,y,\qquad \|\nabla V(x)-\nabla V(y)\| \leq L \|x-y\|,
\]
or, in optimization terminology, where the objective function is \emph{convex and \(L\)-smooth}. For our
purposes here, we shall say that an interval \((0,h_c]\), \(h_c=h_c(L)\), is an \emph{interval of convex
contractivity} of a given RK scheme if, for \(h\in(0,h_c]\), any \(L\)-smooth convex \(V\), and any two
initial points \(\xnot\), \(\xtnot\), the corresponding RK solutions after one time step satisfy
\begin{equation}\label{eq:aim}
 \|\xtone-\xone\|\leq \|\xtnot-\xnot\|.
\end{equation}
By analogy with the result by Humphries and Stuart quoted above, one may perhaps expect that each (consistent)
Runge-Kutta method would possess an interval of convex contractivity; however this is not true. We establish in
Section 3 that the  familiar second-order method due to Runge that for the general system \eqref{eq:generalode}
takes the form
\begin{equation}\label{eq:runge}
 y_1 = y_0 + hF\Big(y_0+\frac{h}{2}F(y_0)\Big)
\end{equation}
possesses no interval of convex contractivity. The proof  proceeds in two stages. We first follow the approach in
\cite{LRP,FRMP}, based on ideas from robust control theory, and identify, for given \(h\) and \(L\), initial
points \(\xnot\), \(\xtnot\) and gradient values
 \[
 \nabla V(\xnot),\quad \nabla V(\xtnot),\quad \nabla
V\Big(\xnot-\frac{h}{2}\nabla V(\xnot)\Big),\quad \nabla V\Big(\xtnot-\frac{h}{2}\nabla V(\xtnot)\Big)
\]
that ensure that \eqref{eq:aim} is violated. In the second stage we provide a counterexample by constructing a
suitable \(L\)-smooth \(V\) by convex interpolation; this is not an easy  task because multivariate convex
interpolation problems with scattered data are difficult to handle \cite{C,CF}. In order not to stop the flow of the paper, some proofs and technical details have been postponed to the final Sections 4
and 5.

For general systems \eqref{eq:generalode}, Dahlquist and Jeltsch \cite{DahlJelt}  considered in an unpublished report 
(summarized in \cite[Chapter 6]{DK}) the monotonicity requirement
\begin{equation}\label{eq:dj}
\forall x,y,\qquad \langle F(x)-F(y), x-y\rangle \leq -\alpha \| F(x)-F(y)\|^2,
\end{equation}
that should be compared with \eqref{eq:nu}. Under this requirement,
they provided \modkzz{a characterization sufficient and necessary
condition) for contractivity of non-confluent Runge--Kutta methods in the setting
of equations $\dot{x}=F(t,x)$ satisfying the monotonicity condition \eqref{eq:dj}.}
Since  it is well known \cite[Theorem 2.1.5]{N14} that \(V\) is \emph{convex and \(L\)-smooth} if and only if
\begin{equation}\label{eq:constraint}
 \forall x,y,\qquad \frac{1}{L}  \| \nabla V(x)-\nabla V(y)\|^2 \leq
\langle \nabla V(x)-\nabla V(y), x-y\rangle,
\end{equation}
it turns out that
 convex, L-smooth gradient systems \eqref{eq:gradient} satisfy \eqref{eq:dj} with \(\alpha = 1/L\) and the Dahlquist-Jeltsch result may be used to derive sufficient conditions for contractivity in our context; in particular
it is possible for some explicit RK schemes to have nonempty intervals of convex contractivity. \modkz{Similar time-step restrictions for explicit RK methods  appear when instead of contractivity \modkzz{one is seeking} to preserve monotonicity \cite{H05}}. For completeness
 we present in Section 2  a version of the theorem by Dahlquist and Jeltsch tailored to our setting of \(L\)-smooth gradient systems. Dahlquist and Jeltsch also proved an opitimality property of Euler's rule among explicit methods and we provide a new proof of their result. \modkz{Optimality of methods of higher order was studied in \cite{K91}.}

Before closing the introduction we point out that there has been much recent interest
\cite{ERRS18,SRB17,SBC16,WWJ16} in interpreting optimization algorithms as discretizations of differential
equations (not necessarily of the form \eqref{eq:gradient}), among other things because differential equations
help to gain intuition on the behaviour of discrete algorithms.

\section{Sufficient conditions for contractivity}

 The application of the \(s\)-stage RK method with coefficients \(a_{ij}\) and weights \(b_j\), \(i,j =
1,\dots,s\), to the system of differential equations \eqref{eq:gradient} results in the relations
\begin{eqnarray}\label{eq:rk}
x_1 &=& x_0+h \sum_{j=1}^s b_j k_j,\\\nonumber
X_i &=& x_0 + h \sum_{j=1}^s a_{ij} k_j,\quad i = 1,\dots,s,\\
k_j &=& -\nabla V(X_j), \quad j = 1,\dots,s. \nonumber
\end{eqnarray}
Here the \(X_i\) and \(k_i\) are the stage vectors and slopes respectively. Of course, the scheme is
consistent/convergent provided that \(\sum_j b_j = 1\).

Item 1 in the Theorem below is essentially Theorem 4.1 in \cite{DahlJelt} and holds for general systems \eqref{eq:generalode} that satisfy \eqref{eq:dj} with \(\alpha = 1/L\) (in fact the proof presented below applies to that more general setting).
The \(s\times s\) symmetric matrix with entries
\[
m_{ij} =b_ia_{ij}+b_ja_{ji}-b_ib_j
\]
that appears in the hypotheses plays a central role in the study of algebraic stability as defined by
Burrage and Butcher, \cite[Definition 12.5]{HaW96} or  \cite[Definition 357B]{B16} and also in symplectic
integration \cite{SS}.

\begin{theorem}
\label{th:main} Let the scheme \eqref{eq:rk} be applied to the gradient system \eqref{eq:gradient} with
convex, \(L\)-smooth \(V\).

Assume that:
\begin{enumerate}
\item The weights \(b_j\), \(j=1,\dots,s\), are nonnegative.
\item The \(s\times s\) symmetric matrix \(\overline M(h)\) with entries \[\overline
    m_{ij}(h)=\frac{2hb_i}{L}\delta_{ij}+h^2m_{ij}\] (\(\delta\) is Kronecker's symbol) is positive
    semidefinite.
\end{enumerate}
Then:
\begin{enumerate}
 \item If \(x_1\) and \(\widetilde x_1\) are the RK solutions after a step of lenght \(h>0\) starting from
\(\xnot\) and \(\xtnot\) respectively the contractivity estimate \eqref{eq:aim} holds.

\item In particular, if
\(x^\star\) is a minimizer of \(V\), then
\[
\|x_1-x^\star\|\leq \|x_0-x^\star\|.
\]
\end{enumerate}
\end{theorem}

\emph{Proof.} We start from the identity \cite[Theorem 12.4]{HaW96}
\[
\| \widetilde x_1-x_1\|^2 = \| \widetilde x_0-x_0\|^2+2h \sum_{i=1}^s b_i \langle \widetilde k_i-k_i, \widetilde X_i- X_i\rangle- h^2 \sum_{i,j=1}^sm_{ij} \langle \widetilde k_i-k_i,\widetilde k_j-k_j \rangle,
\]
where \(\widetilde X_i\) and \(\widetilde k_i\) respectively denote the stage vectors and slopes for the step
\(\widetilde x_0\mapsto \widetilde x_1\). (This identity holds if \(\langle \cdot,\cdot\rangle\) and
\(\|\cdot\|\) are replaced by any symmetric bilinear map and the associated quadratic map respectively, see
\cite[Lemma 2.5]{SS}.) From \eqref{eq:constraint}, for \(i= 1,\dots, s\),
\[
\langle \widetilde k_i-k_i, \widetilde X_i- X_i\rangle \leq -\frac{1}{L}
\langle \widetilde k_i-k_i, \widetilde k_i- k_i\rangle,
\]
which implies, in view of the nonnegativity of the weights,
\modk{\[
\| \widetilde x_1-x_1\|^2 \leq \| \widetilde x_0-x_0\|^2 -\sum_{i,j=1}^s \overline m_{ij}(h)
\langle \widetilde k_i-k_i,\widetilde k_j-k_j \rangle.
\]}
If \(\overline M(h)\) is positive semidefinite the sum is \(\geq 0\) and the proof is complete. In addition, if we now set $\widetilde{x}_{0}=x^{\star}$, we trivially obtain $\|x_1-x^\star\|\leq \|x_0-x^\star\|$.
  \bigskip

We next present some examples; the interested reader may find a full discussion in the report \cite{DahlJelt}. Hereafter \(Q\succeq 0\) means that the matrix \(Q\) is positive semidefinite.

\emph{Example 1.} For Euler's rule, \(s=1\), \(a_{11} =0\), \(b_1=1\), we find \(\overline M(h) = 2h/L-h^2\)
and therefore we have contractivity for \(h\) in the interval \((0,2/L]\). This happens to coincide with the
familiar stability interval for the linear scalar test equation \((d/dt) x = -Lx\), \(L>0\). The restriction
\(h\leq 2/L\) on the step size  is well known in the analysis of the gradient descent algorithm, see e.g.\
\cite{N14}. Observe that the scalar test equation arises from the \(L\)-smooth convex potential \(V= Lx^2/2\)
and that therefore no RK scheme can have an interval of convex contractivity longer than its linear stability
interval.

\emph{Example 2.} The formula two-stage, second order \eqref{eq:runge} presented in the introduction has
\(b_1=0\), \(b_2 = 1\) and \(a_{21} = 1/2\). Thus
\[\overline M(h) =
\left[
\begin{matrix}
0 & \frac{h^2}{2}\\
\frac{h^2}{2} & \frac{2h}{L}-h^2
\end{matrix}
\right].
\]
There is  no value of \(h>0\) for which this matrix is \(\succeq 0\). In
Theorem~\ref{th:convex} we shall show that the scheme has no interval of convex contractivity. Hence for this RK method the sufficient condition in Theorem~\ref{th:main} is actually \emph{necessary}. Note the necessity, under the requirement \eqref{eq:dj}, of the hypotheses of Theorem 4.1 in \cite{DahlJelt} was not discussed by Dahlquist and Jeltsch.

\emph{Example 3.} Explicit, two-stage, first-order scheme with \(b_1=b_2 = 1/2\) and \(a_{21} = 1/2\).   Here
\[\overline M(h)=
\left[
\begin{matrix}
\frac{h}{L}-\frac{h^2}{4} & 0\\
0 & \frac{h}{L}-\frac{h^2}{4}
\end{matrix},
\right]
\]
 and we have contractivity for \(0<h\leq 4/L\). This could have been concluded from  Example 1, because
 performing one step with this method yields the same result as taking  two steps of length \(h/2\) with Euler's rule and accordingly, for this method,
  \(h/2 \leq 2/L\) ensures contractivity.

\emph{Example 4.} We may generalize as follows. Consider the explicit \(s\)-stage, first-order scheme with
Butcher tableau
\begin{equation}\label{eq:euler}
\begin{matrix}
0 & 0 & 0 & \cdots&0\\
b_1 & 0 & 0& \cdots&0\\
b_1 & b_2 & 0& \cdots&0\\
\vdots &\vdots&\vdots& \vdots&\vdots\\
b_1 & b_2 & b_3& \cdots&0\\\hline
b_1 & b_2 & b_3& \cdots&b_s
\end{matrix}
\end{equation}
(i.e., \(a_{ij} = b_j\) whenever \(i>j\)) with
\[ \sum_{i= 1}^sb_i= 1,\qquad b_i \geq 0,\quad i = 1,\dots,s.\]
 Performing one step with this scheme is equivalent to successively performing \(s\) steps with Euler's rule with step-sizes \(b_1h\), \dots , \(b_sh\), and therefore contractivity is ensured \modkz{in the case when} \(h\max_i b_i\leq 2/L\). This conclusion may alternatively be reached by applying \modkz{Theorem \ref{th:main}}; the method has \(\overline M(h)\) given by
\begin{equation}\label{eq:M}
 {\rm diag}\big( 2hb_1/L-h^2b_1^2,\quad 2hb_2/L-h^2b_2^2,\: \dots,\:  2hb_s/L-h^2b_s^2\big),
\end{equation}
a matrix that is \(\succeq 0\) if and only if \(h\max_i b_i\leq 2/L\). If we see the weights as parameters,
then the least severe restriction on \(h\) is attained by choosing equal weights \(b_i = 1/s\), \(i=1,\dots,
s\), leading to the condition \(h\leq 2s/L\). But then one is really time-stepping with Euler rule with
stepsize \(h/s\).\bigskip

Recall that RK schemes are called reducible if they give the same numerical results as a scheme with fewer stages; reducible methods are then completely devoid of interest.
It is not difficult to prove (see \cite[Corollary 3.4]{DahlJelt}) that RK schemes that are not reducible and for which \(\overline M(h)\succeq 0\) for at least one value of \(h\) have all its weights strictly positive. It is also known that
irreuducible, explicit methods with positive weights have order \(\leq 4\), \cite[Theorem 4.4]{DahlJelt}.

The next result is essentially Theorem 5.1 in \cite{DahlJelt} and shows that among explicit methods Euler's rule has the longest interval of convex contractivity if intervals are scaled in terms of the number of stages so as to take the amount of work per step. Our purely algebraic proof is different from the analytic one given by Dahlquist and Jeltsch. Note that, in view of the comment we just made, the weights are assumed to be \(>0\).

\begin{theorem}\label{th:optimal}
  Consider an \(s\)-stage, explicit, consistent RK method with  weights \(>0\).
\begin{enumerate}
 \item If for some \(h>0\),  \(\overline M(h)\succeq 0\), then \(h\leq 2s/L\).
 \item If for \(h=2s/L\), \(\overline  M(h)\succeq 0\), then the method is necessarily given by
     \eqref{eq:euler} with \(b_i = 1/s\), \(i=1,\dots,s\) (i.e., it is the concatenation of \(s\) Euler
     substeps of equal length \(h/s\)).
 \end{enumerate}
 \end{theorem}

 \emph{Proof.} For the first item, we first note that, as we saw in Example 4, the result is true for the particular case where the scheme
 is of the  form \eqref{eq:euler}, i.e.,  a concatenation of Euler's substeps. Let
 \(\overline M_\star(h)\) be the matrix associated with the scheme of the form \eqref{eq:euler} that possesses the same weights as the given scheme (recall that this matrix was computed in \eqref{eq:M}). The first item will be proved if we show that \(\overline M(h)\succeq 0\) implies
 \(\overline M_*(h)\succeq 0\), because, as we have just noted, the last condition guarantees that \(h\leq 2s/L\).  Assume  that \(\overline M(h)\succeq 0\). Then, its diagonal entries must be nonnegative,
 \[0\leq \overline m_{ii}(h) = 2hb_i/L-h^2b_i^2,\qquad i=1,\dots, s,\]
 and, in view of \eqref{eq:M}, this entails that \(\overline M_\star(h)\succeq 0\), as we wanted to establish.

 We now prove the second part of the theorem. If \(\overline M(2s/L)\succeq 0\), then
 \[
 0\leq \overline m_{ii}(2s/L) = 4sb_i/L^2-4s^2b_i^2/L^2,\qquad i=1,\dots, s,
 \]
 or, after dividing by  \(4b_is^2/L^2>0\), \(b_i\leq 1/s\). Since \(\sum_{i=1}^s b_i = 1\), we conclude that
 \(b_i =1/s\), \(i=1,\dots,s\), which leads to  \(\overline m_{ii}(2s/L)=0\) for each \(i\). A  semidefinite positive matrix
 with vanishing diagonal elements must be the null matrix and therefore for \(i>j\)
 \[ 0 = \overline m_{ij}(2s/L) = (2s/L)^2 (b_ia_{ij}-b_ib_j)\]
 and then \( a_{ij} = b_j\). The proof is now complete.
\bigskip

\section{An RK scheme without convex contractivity interval}
\label{s:control}
In this section we show that the RK scheme \eqref{eq:runge} has no interval of convex contractivity.

 For the system \eqref{eq:gradient}, we write the formulas for performing one step from the initial
points \(\xnot\) and \(\xtnot\) in \(\R^d\) as
\begin{equation}\label{eq:step}
\xone = \xnot + h \kh,\quad \xtone = \xtnot + h \kth,
\quad
\xh = \xnot +\modkzz{ \frac{h}{2}}  \knot,\quad \xth = \xtnot + \frac{h}{2} \ktnot,
\end{equation}
with
\begin{equation}\label{eq:feedback}
\knot = -\nabla V(\xnot),\quad
\ktnot = -\nabla V(\xtnot),\quad
\kh = -\nabla V(\xh),\quad
\kth = -\nabla V(\xth)
\end{equation}
(the subindices \(0\), \(1\), \(h\) refer to the beginning of the step, \(t=0\), the end of the step, \(t=h\),
and the halfway location, \(t=h/2\), respectively). Following the approach in \cite{LRP,FRMP}, we regard
\(\xnot\), \(\xtnot\), \(\knot\), \(\ktnot\), \(\kh\), \(\kth\), as \emph{inputs}, and \(\xnot\), \(\xtnot\),
\(\xh\), \(\xh\), \(\xone\), \(\xtone\) as \emph{outputs}\footnote{\modk{Note that \(\xnot\), \(\xtnot\) are both inputs and outputs.}}. The relations \eqref{eq:feedback} provide a
\emph{feedback} loop that expresses the inputs \(\knot\), \(\ktnot\), \(\kh\), \(\kth\) as values of a
nonlinear  function \(\phi=-\nabla V\) computed at  the outputs \(\xnot\), \(\xtnot\), \(\xh\), \(\xth\). The
function \(\phi\) that establishes this feedback is  the negative gradient of some \(V\) that  is convex and
\(L\)-smooth. According to \eqref{eq:constraint},
this implies that the  vectors \(\knot\), \(\ktnot\), \(\kh\), \(\kth\) delivered by the feedback loop must
obey the following constraints:
\begin{eqnarray}
\label{eq:contraint 1}
\frac{1}{L} \| \ktnot-\knot\|^2 & \leq & -\langle \ktnot-\knot, \xtnot-\xnot\rangle,\\
\label{eq:contraint 2}
\frac{1}{L} \| \kth-\kh\|^2 & \leq & -\langle \kth-\kh, \xth-\xh\rangle,\\
\label{eq:contraint 3}
\frac{1}{L} \| \kh-\knot\|^2 & \leq & -\langle \kh-\knot, \xh-\xnot\rangle,\\
\label{eq:contraint 4}
\frac{1}{L} \| \kth-\ktnot\|^2 & \leq & -\langle \kth-\ktnot, \xth-\xtnot\rangle,\\
\label{eq:contraint 5}
\frac{1}{L} \| \kth-\knot\|^2 & \leq & -\langle \kth-\knot, \xth-\xnot\rangle,\\
\label{eq:contraint 6}
\frac{1}{L} \| \kh-\ktnot\|^2 & \leq & -\langle \kh-\ktnot, \xh-\xtnot\rangle
\end{eqnarray}
(we are dealing with four gradient values and therefore \eqref{eq:constraint} may be applied in \({4\choose
2}=6\) ways). In a \emph{robust} control approach, we will not assume at this stage that the vectors \(k\) are
values of one and the same function \(-\nabla V\), on the contrary the vectors \(k\) are seen as arbitrary
except for the above constraints. More precisely, for fixed \(L\) and \(h\), we investigate the lack of
contractivity by studying the real function
\begin{equation}\label{eq:objective}
\frac{\|\xtone-\xone\|^2}{\|\xtnot-\xnot\|^2}
\end{equation}
of the input variables \(\xnot\), \(\xtnot\), \(\knot\), \(\ktnot\), \(\kh\), \(\kth\), subject to the
constraints \(\xtnot\neq \xnot\) and \eqref{eq:contraint 1}--\eqref{eq:contraint 6}. Here \(\xh\), \(\xh\),
\(\xone\), \(\xtone\) are known  linear combinations of the inputs given in \eqref{eq:step}.

Our task is made easier by the following observations. First of all, multiplication of  \(\xnot\), \(\xtnot\),
\(\xh\), \(\xth\), \(\xone\), \(\xtone\), \(\knot\), \(\ktnot\), \(\kh\), \(\kth\) by the same scalar
\(\lambda
>0\) preserves the relations \eqref{eq:step}, the constraints \eqref{eq:contraint 1}--\eqref{eq:contraint 6}
and the value of the quotient \eqref{eq:objective}. Therefore we may assume at the outset that
\(\|\xtnot-\xnot\|= 1\). In addition, since the problem is also invariant by translations and rotations in
\(\R^d\), we may set \(\xnot = 0\in\R^d\) and \(\xtnot = [1,0,0,\dots,0]^T\). After these simplifications, we
are left with the task of ascertaining if we can make  \(\|\xtone-\xone\|^2\) larger than \(1\) by choosing
appropriately the vectors \(\knot\), \(\ktnot\), \(\kh\), \(\kth\) subject to the constraints. Here is a
choice in \(\R^2\) that works (see Section 4 for the origin of these vectors)
\begin{eqnarray}
\label{eq:magic1}
\knot &=& [0, -3/h]^T,\\
\label{eq:magic2}
\ktnot &=& [-L/2, -3/h+L/2]^T,\\
\label{eq:magic3}
\kh&=& [0, -3/h+L]^T,\quad\\
\label{eq:magic4}
\kth &=& [L^3h^2/64, -3/h+L -L^2h/8]^T.
\end{eqnarray}
In fact, with
\begin{equation}\label{eq:future0}
\xnot = [0,0]^T,\qquad \xtnot = [1,0]^T
\end{equation}
and \eqref{eq:magic1}--\eqref{eq:magic4}, the relations \eqref{eq:step} yield
\begin{eqnarray}\label{eq:future1}
 \xh &=& [0,-3/2]^T,\\
 \label{eq:future2}
 \xth &=& [1-Lh/4,-3/2+Lh/4]^T,\\
 \label{eq:future3}
 \xone &=& [0,-3+Lh]^T,
 \\
 \label{eq:future4}
 \xtone &=& [1+L^3h^3/64,-3+Lh+L^2h^2/8]^T.
\end{eqnarray}
It is a simple exercise to check that the constraints are satisfied
 at least for \(Lh\leq 3\). In addition
 \[ \xtone-\xone = [1+L^3h^3/64, L^2h^2/8]^T
 \]
 and, accordingly,
\begin{equation}\label{eq:growth}
\|\xtone-\xone\|^2 = 1+\frac{1}{32} L^3h^3+ \frac{1}{64}L^4h^4+\frac{1}{4096}L^6h^6>1= \|\xtnot-\xnot\|^2.
\end{equation}
(The third power in \(h^3\) matches the size of the local error of the scheme.)

\modk{\begin{remark}
The vectors \eqref{eq:magic1}--\eqref{eq:magic4} become longer as \(h\) decreases. This is a consequence of
the way we addressed the study of \eqref{eq:objective} where we fixed the length of \(\xtnot-\xnot\) for
mathematical convenience. As pointed out above we could alternatively have chosen \(\xnot =[0,0]^T\),
\(\xtnot=[h,0]^T\) and multiplied \eqref{eq:magic1}--\eqref{eq:magic4} by a factor of \(h\) and that would
have given a configuration with bounded gradients resulting in lack of contractivity.
\end{remark}}

\begin{figure}
\begin{center}
\vspace{-4.5cm}
\includegraphics[scale=0.5]{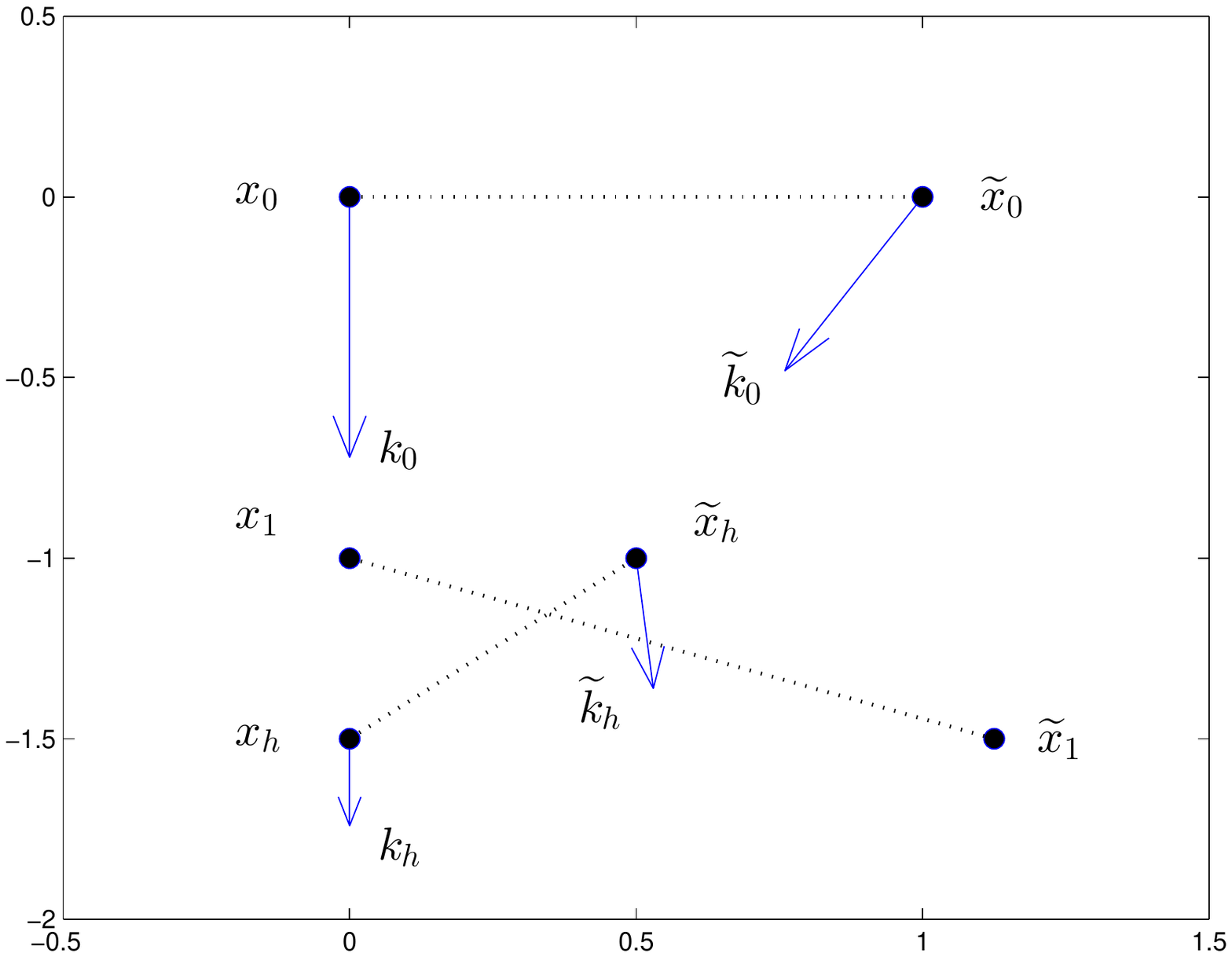}
\end{center}
\vspace{-3cm}
\caption{A configuration that satisfies the constraints resulting from convexity and \(L\)-smoothness and leads to lack of contractivity for \(L=2\), \(h=1\)}
\label{fig:insight}\end{figure}

To get some insight,  we have depicted in Figure~\ref{fig:insight}, when \(L=2\), \(h=1\), the points
\(\xnot\), \(\xtnot\), \(\xh\), \(\xth\), \(\xone\), \(\xtone\)  along with the vectors \(\knot\), \(\ktnot\),
\(\kh\), \(\kth\) (for clarity, the vectors have been drawn after multiplying their length by \(0.8\)). The
difference vector \(\ktnot-\knot\) forms, as required by convexity, an \emph{obtuse} angle with
\(\xtnot-\xnot\) and this causes \(\xth-\xh= \xtnot-\xnot + (h/2) (\ktnot-\knot)\) to be \emph{shorter} than
\(\xtnot-\xnot\). Similarly the difference \(\kth-\kh\) forms by convexity an \emph{obtuse} angle with
\(\xth-\xh\) and if \(\xone\) and \(\xtone\) were alternatively defined as \(\xh+(h/2)\kh\) and
\(\xth+(h/2)\kth\) respectively
 we see from the Figure that we would have \(\|\xtone-\xone\|\leq  \|\xtnot-\xnot\|\).
(That alternative time stepping was studied in Example 3 in the preceding section.) However for our RK scheme
\eqref{eq:runge} the direction of \(\kh\) is used to displace \(\xnot\) (rather than \(\xh\)) to get \(\xone\)
(and similarly for the points with tilde); the vector \(\kth-\kh\) forms an \emph{acute} angle with
\(\xtnot-\xnot\) and this makes it possible for \(\xtone-\xone\) to be longer than \(\xtnot-\xnot\). For
smaller values of \(L\) and/or \(h\) the effect is not so marked as that displayed in the figure but is
nevertheless present.

While \eqref{eq:growth} is consistent with the scheme having no interval of convex contractivity, we are not yet
done, because it is not obvious whether there is a convex, \(L\)-smooth \(V\) that realizes the relations
\eqref{eq:feedback} for the \(x\)'s and \(k\)'s we have found. Nevertheless the preceding material will
provide the basis for proving in the final section the following result:
\begin{theorem}\label{th:convex}
Fix \(L>0\). For the RK scheme \eqref{eq:runge} and each arbitrarily small value of \(h>0\), there exist an
\(L\)-smooth, convex \(V\) and initial points \(\xnot\) and \(\xtnot\) such that \eqref{eq:aim} is not
satisfied. As a consequence the scheme does not possess an interval of convex contractivity.
\end{theorem}
One could perhaps say that the method has an \emph{empty} interval of convex contractivity.

\section{The construction of the auxiliary gradients}

The proof of Theorem~\ref{th:convex} hinges on the use  of the vectors \modk{\eqref{eq:magic1}--\eqref{eq:magic4}}. In this section we briefly describe how we constructed them.

Let us introduce the vectors in \(\R^2\)
\[
\dnot = \xtnot-\xnot,\qquad \deltah = \xth-\xh, \qquad \done = \xtone-\xone
\]
and
\[
\Dnot = \ktnot-\knot, \qquad \Dh =\kth-\kh,
\]
so that \(\done = \dnot+h\Dh\) and \(\delta_h = \dnot+(h/2)\Dnot\).
We fixed \(\dnot=[1,0]^T\) as explained in Section \ref{s:control}, saw \(\Dnot\) and \(\Dh\) as
variables in \(\R^2\) and considered the problem of maximizing \(\|\done\|^2\) under the constraints
\eqref{eq:contraint 1}--\eqref{eq:contraint 2}, i.e
\[
\frac{1}{L} \|\Dnot\|^2 \leq -\langle \Dnot,\dnot\rangle,\qquad
\frac{1}{L} \|\Dh\|^2 \leq -\langle \Dh,\deltah\rangle,
\]

With some patience, we solved this maximization problem  analytically in closed form after
 introducing Lagrange multipliers.  Both
constraints are active at the solution. The expression of the maximizer  is a complicated function
 of \(L\) and \(h\) and to simplify the subsequent algebra we expanded that expression in powers
 of \(h\) and kept the leading terms. This resulted in
\[
\Dnot = [-L/2, L/2]^T, \qquad \Dh = [L^3h^2/64, - L^2h/8]^T
\]
(there is a second solution obtained by reflecting this with respect to the first coordinate axis).

Once we had found candidates for the differences \(\ktnot-\knot\), \(\kth-\kh\), we identified suitable
candidates for \(\knot\) and \(\kh\). We arbitrarily fixed the direction of \(\knot\) by choosing it to be
perpendicular to \(\delta_0\) (see \eqref{eq:magic1}). Its second component was sought in the form \(c/h\)
(\(c\) a constant) so that the distance between \(\xh\) and \(\xnot\) behaved like \(\mathcal{O}(1)\) as
\(h\downarrow 1\) (recall that we have scaled things in such a way that \(\xtnot\) and \(\xnot\) are also at a
distance \(\mathcal{O}(1)\) as \(h\downarrow 1\)).  We also took \(\kh\) to be perpendicular to \(\delta_0\);
the second component of this vector was chosen to be of the form \(c/h-c^\prime L\) so as to have \(\kh-\knot=
-c^\prime L\) with a view to satisfying \eqref{eq:contraint 3}. After some numerical experimentation we saw that the
values \(c=3\), \(c^\prime = 1\) led to a set of vectors for which all six contraints \eqref{eq:contraint
1}--\eqref{eq:contraint 6} hold at least for \(Lh\leq 3\).

For the sake of curiosity we also carried out numerically the maximization of \eqref{eq:objective} subject to
the constraints. It turns out that the maximum value of the quotient is approximately \(1+0.032 L^3h^3\) for
\(h\) small, independently of the dimension \(d\geq 2\) of the problem (for \(d=1\) the experiments suggest
that the scheme is contractive). Since, in \eqref{eq:growth}, \(1/32 = 0.03125\) the vectors
\eqref{eq:magic1}--\eqref{eq:magic4} are very close to providing the combination of gradients that leads to
the greatest dilation \eqref{eq:objective}.

\section{Proof of Theorem \ref{th:convex}}
The proof proceeds in two stages. We first construct an auxiliary piecewise linear, convex \(\modk{\widetilde{V}}\) and then we
regularize it to obtain \(V\).
\subsection{Constructing a piecewise linear potential by convex interpolation}
Let \(L>0\) be the Lispchitz constant and set \(\Lp = \alpha L\), where \modkzz{$\alpha$} is a positive safety factor,
independent of \(L\) and \(h\), whose value will be determined later. Restrict hereafter the attention to
values of \(h\) with \(h\Lp\leq 1\). We wish to construct a potential \(\modk{\widetilde{V}}\) for which the application of the
RK scheme starting from the two initial conditions \eqref{eq:future0} lead to the relations
\eqref{eq:magic1}--\eqref{eq:magic4}, \eqref{eq:future1}--\eqref{eq:future4} \emph{with \(\Lp\) in lieu of
\(L\)} and therefore, \modk{as we know}, to lack of contractivity.

We consider the following four (pairwise distinct) points in the plane \(\R^2\) of the variable \(\zeta\) (see
\eqref{eq:future1}--\eqref{eq:future4})
\begin{eqnarray*}
Z_1&=& [0,0]^T,\\
Z_2 &=& [1,0]^T,\\
Z_3&=& [0,-3/2]^T,\\
Z_4 &=& [1-\Lp h/4,-3/2+\Lp h/4]^T,
\end{eqnarray*}
and associate with them the four (pairwise distinct) vectors (see \eqref{eq:magic1}--\eqref{eq:magic4})
\begin{eqnarray*}
\modkzz{G_1} &=& [0, 3/h]^T,\\
\modkzz{G_2} &=& [\Lp/2, 3/h-\Lp/2]^T,\\
\modkzz{G_3}&=& [0, 3/h-\Lp]^T,\quad\\
\modkzz{G_4} &=& [-\Lp^3h^2/64, 3/h-\Lp +\Lp^2h/8]^T,
\end{eqnarray*}
and four real numbers \(F_i\) that will be determined below. We then pose the following \emph{Hermite convex
interpolation problem:} Find
 a real convex function \(\modk{\widetilde{V}}\) defined in  \(\R^2\), differentiable in the neighbourhood of the \(Z_i\), and
 such that
 \[
 \modk{\widetilde{V}}(Z_i) = F_i,\qquad \nabla \modk{\widetilde{V}}(Z_i) = G_i,\qquad i = 1,\dots,4.
 \]

 If the interpolation problem has a solution, then the tangent plane to \(\eta=\modk{\widetilde{V}}(\zeta)\) at \(Z_i\) is
 given by the equation
 \( \eta = \pi_i(\zeta)\) with
 \[
 \modk{\pi_i(\zeta)}\ = F_i+\langle G_i, \zeta-Z_i\rangle, \qquad i=1,\dots,4.
 \]
 and, by convexity,
 \begin{equation}\label{ec:necesaria}
 F_i \geq \pi_j(Z_i), \qquad i\neq j,\quad i,j = 1,\dots,4.
 \end{equation}
 This is then a necessary condition for the Hermite problem to have a solution. We found the
 following set of values
 \begin{eqnarray*}
 F_1 &=& 0,\\
 F_2 &=& \frac{\Lp}{4},\\
 F_3 & = &-\frac{9}{2h}+\frac{9\Lp}{8},\\
 F_4 & = & -\frac{9}{2h}+\frac{15\Lp}{8}-\frac{ \Lp^2h}{4}+\frac{\Lp^3h^2}{128},
 \end{eqnarray*}
that  satisfy the relations \eqref{ec:necesaria} (in fact they satisfy all of them with strict inequality).

It is not difficult to see \cite{C,CF}, that once we have ensured \eqref{ec:necesaria},
 the Hermite problem is solvable. The solution is not unique and among all solutions the minimal is clearly
 given by the piecewise linear function
 \[
 \modk{\widetilde{V}}(\zeta) = \max\{ \pi_i(\zeta) : i = 1,\dots, 4\}.
 \]

 From Section~\ref{s:control} we conclude that, if the RK scheme is applied to solve the gradient system
 associated with \(\modk{\widetilde{V}}\) with starting points \(\xnot = Z_1\), \(\xtnot = Z_2\),
 then \eqref{eq:growth} holds with \(L\) replaced by \(\Lp\) and there is no contractivity.
 However, the proof is not complete because \(\modk{\widetilde{V}}\) is not continuously differentiable (let alone \(L\)-smooth).
Accordingly we shall regularize \(\modk{\widetilde{V}}\) to construct the potential \(V\) we need.
\begin{figure}
\begin{center}
\vspace{-4.5cm}
\includegraphics[scale=0.5]{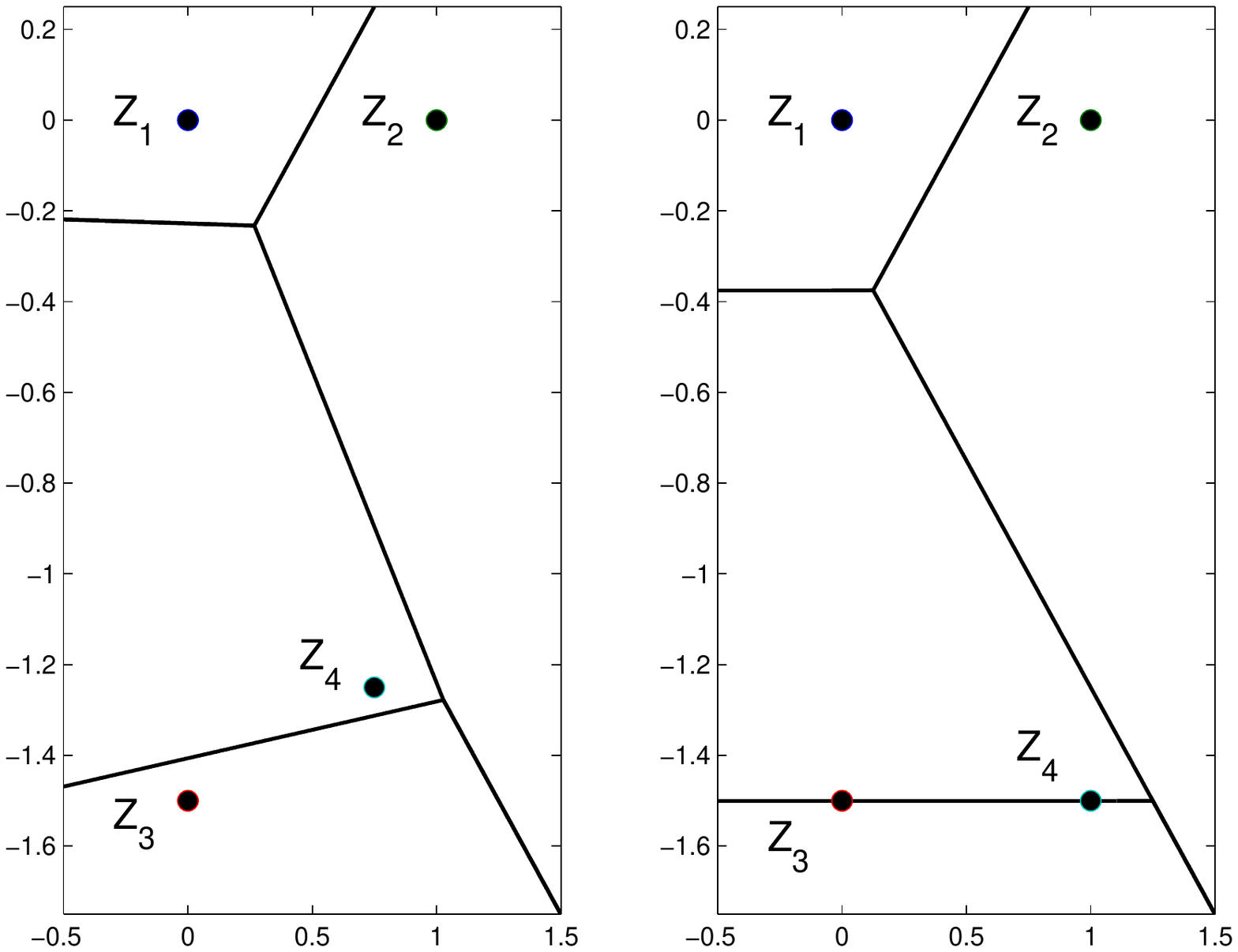}
\end{center}
\vspace{-3cm}
\caption{Left: points \(Z_i\) and tessellation associated with the piecewise linear convex interpolant when
\(\Lp h=1\). Right: points \(Z_i\) and tessellation in the limit \(\Lp h\downarrow 0\)}
\label{fig:tessellation}\end{figure}

Before we do so, it is convenient to notice that \(\modk{\widetilde{V}}\) gives rise to four closed, convex regions \cite{C,CF}
\[
{\cal R}_i = \{ \zeta: \modk{\widetilde{V}}(\zeta) = \pi_i(\zeta)\}, \quad i = 1,\dots, 4,
\]
that tessellate the plane. Clearly \(Z_i\in{\cal R}_i\), \(i=1,\dots,4\). The equations of the lines that
bound the regions are of course found by intersecting the planes \( \eta = \pi_i(\zeta)\), \(i=1,\dots,4\).
After carrying out the corresponding trite computations, it turns out that those boundaries depend on \(h\)
and \(\Lp\) only through the product \(\Lp h\). (By the way, the same is true of the coordinates of the points
\(Z_i\).) For \(\Lp h=1\), the maximum value under consideration of the product \(\Lp h\), we have depicted
the interpolation nodes and regions in the left panel of Figure~\ref{fig:tessellation}. Note that the gradient
\(\nabla \modk{\widetilde{V}}\) takes the constant value \(G_i\) in the interior of the region \({\cal R}_i\). This gradient is
then discontinuous at the boundaries of the tessellation; from the analytic expressions for the \(G_i\) we see
that the jumps \(\|G_i-G_j\|\) at the boundaries may be bounded above by \(C_1\Lp\) with \(C_1\) a constant
independent of \(\Lp\) and \(h\).

While the interpolation problem above only makes sense for positive \(h\), the points \(Z_i\) and the
tessellation have well-defined limits as \(h\downarrow 0\); these limits are depicted in the right panel of
Figure~\ref{fig:tessellation}.
 Note for future reference that, in the limit, \(Z_3\) and \(Z_4\) are on the common boundary of
\({\cal R}_3\) and \({\cal R}_4\).
\subsection{Regularization by convolution}
For \(\zeta\in\R^2\) let us denote by \({\cal S}(\zeta)\subset\R^2\) the closed square centered at \(\zeta\)
with side \(\ell/2\) (i.e. the closed \(L_\infty\)-ball centered at  \(\zeta\) with radius \(\ell/2\))). The
regularization procedure uses the real-valued function \(\chi(\zeta)\) such that \(\chi(\zeta) = 1/\ell^2\) if
\(\zeta\in{\cal S}(0)\) and \(\chi(\zeta) = 0\) if \(\zeta\notin{\cal S}(0)\). Clearly \(\int_{\R^2}
\chi(\zeta)\, d\zeta=1\).

We fix the value of \(\ell\) in such a way that for all \(\Lp>0\) and all \(h\leq 1/\Lp\) (see
Figure~\ref{fig:tessellation})
\[
{\cal S}(Z_1)\subset{\cal R}_1,\quad
{\cal S}(Z_2)\subset {\cal R}_2 ,\quad
{\cal S}(Z_3)\subset {\cal R}_3 \cup {\cal R}_4,\quad
{\cal S}(Z_4)\subset {\cal R}_3 \cup {\cal R}_4;
\]
it is not possible to achieve \({\cal S}(Z_3)\subset {\cal R}_3  \), or \({\cal S}(Z_4)\subset {\cal R}_4 \)
because \(\ell\) is not allowed to depend on \(h\) and, as \(h\) decreases, \(Z_3\) and \(Z_4\) approach the
boundary of \({\cal R}_3\) and \({\cal R}_4\), as we just pointed out.

We define the regularized potential by the convolution
\[
V(\zeta) = \int_{\R^2} \chi(\zeta^\prime)\, \modk{\widetilde{V}}(\zeta-\zeta^\prime)\,d\zeta^\prime.
\]
Each translated function \(\zeta\mapsto \modk{\widetilde{V}}(\zeta-\zeta^\prime)\) is convex and \(\chi(\zeta^\prime)\geq 0\)
so that \(V\) is convex, as a convex combination of convex functions. Furthermore
\[
\nabla V(\zeta) = \int_{\R^2} \chi(\zeta^\prime)\,\nabla \modk{\widetilde{V}}(\zeta-\zeta^\prime)\,d\zeta^\prime
\]
(the integrand is not defined on the lines that define the tessellation) or
\[
\nabla V(\zeta) = \int_{\R^2} \chi(\zeta-\zeta^\prime)\,\nabla \modk{\widetilde{V}}(\zeta^\prime)\,d\zeta^\prime=
\frac{1}{\ell^2} \int_{ \{\zeta^\prime \in{\cal S}(\zeta)\} } \nabla \modk{\widetilde{V}}(\zeta^\prime)\,d\zeta^\prime.
\]
Since \(\zeta^\prime\mapsto \nabla \modk{\widetilde{V}}(\zeta^\prime)\) is piecewise constant with value \(G_i\) in the
interior of \({\cal R}_i\), \(i=1,\dots,4\), for each fixed \(\zeta\), the vector \(\nabla V(\zeta)\) is a
convex linear combination of the vectors \(G_i\), \(i=1,\dots,4\), and the weights of this combination are
given by \((1/\ell^2)\) times the areas of the intersections \({\cal S}(\zeta)\cap {\cal R}_i\). This shows
that \(\nabla V\) is a continuous function (i.e. that \(V\) is continuous differentiable). In addition,  if
for a given location \(\zeta\) the square \({\cal S}(\zeta)\) is entirely contained in one of the regions
\({\cal R}_{i_0}\), then \(\nabla V(\zeta)=G_{i_0}\). By our choice of \(\ell\) it follows that
\begin{equation}\label{eqz1g1z2g2}
\nabla V(Z_1) = G_1, \qquad \nabla V(Z_2) = G_2.
\end{equation}

The geometric interpretation of the definition of \(\nabla V(\zeta)\) also shows that \(\nabla V\) is
Lipschitz continuous with a Lipschitz constant of the form \(C_2 D/\ell\), where \(D\) is an upper bound for
the size of the jumps \(\| G_i-G_j\|\), \(i,j=1,\dots,4\). As remarked earlier, \(D= C_1\Lp\), so that
\(\nabla V\) is is Lipschitz continuous with Lipschitz constant \( C_1C_2\Lp/\ell\). Therefore by choosing our
safety factor as \(\alpha = \ell/(C_1C_2)\), the potential \(V\) will be convex and \(L\)-smooth.

Finally take RK solutions for the problem \eqref{eq:gradient} from the poins \(\xnot = Z_1\) and \(\xtnot
=Z_2\). From \eqref{eqz1g1z2g2} and the definition of \(G_1\) and \(G_2\), we have \(\xh = Z_3\) and \(\xth = Z_4\). Next
\begin{eqnarray*}
&&\nabla V(\xh) = \nabla V(Z_3) = \lambda G_3+ (1-\lambda) G_4,\\ &&\nabla V(\xth)=\nabla V(Z_4) = (1-\mu) G_3+ \mu G_4
\end{eqnarray*}
where \(\lambda\) is \(1/\ell^2\) times the area of \({\cal S}(Z_3)\cap {\cal R}_3\) and \(\mu\) is
\(1/\ell^2\) times the area of \({\cal S}(Z_4)\cap {\cal R}_4\).  We observe that  \( \lambda > 1/2\) \modk{for $h>0$} because
\({\cal S}(Z_3)\cup {\cal R}_3\) clearly has
 more area than \({\cal S}(Z_3)\cup {\cal R}_4\). Similarly \(\mu > 1/2  \) \modk{for $h>0$}. The quantities \(\lambda\) and
 \(\mu\) depend on \(\Lp\) and \(h\) and approach \(1/2\) as \(h\downarrow 0\).
We then find
\[
\xtone -\xone = [ 1 +  \nu \Lp^3h^3/64, - \nu \Lp^2h^2/8]^T, \qquad \nu = \mu -(1-\lambda)>0
\]
and
\[
\|\xtone -\xone\|^2 = 1+ \frac{1}{32}\nu \Lp^3h^3+ \frac{1}{64}\nu^2 \Lp^4h^4+\frac{1}{4096}\nu^2 \Lp^6h^6 >1.
\]
This estimate is worse  than \eqref{eq:growth} due to the presence of \(\Lp\) and \(\nu\), but still
sufficient to prove the theorem. By using functions \(\chi\) smoother than the one we used above, it is possible to construct by convolution smoother
potentials \(V\). However, our choice here results in a clearer proof.

\bigskip {\bf Acknowledgements.} J.M.S. was supported by project
MTM2016-77660-P(AEI/ FEDER, UE)
 funded by MINECO (Spain). He would like to thank the Isaac Newton Institute for Mathematical
 Sciences for support and hospitality during the programme
 \lq\lq Geometry, compatibility and structure preservation in computational differential equations\rq\rq\
 when work on this paper was undertaken. This work
 was supported by EPSRC Grant Number EP/R014604/1. K.C.Z was supported by the Alan
Turing Institute under the EPSRC grant EP/N510129/1. The authors are also thankful to J. Carnicer (Zaragoza) for bringing to their
attention a number of helpful references on convex interpolation.

\end{document}